%
%

\documentstyle{amsppt}
\pagewidth{160mm}
\pageheight{240mm}
\magnification=1200
\NoBlackBoxes
\def\J{{\Cal J}}
\def\p{\tilde{p}}
\def\q{\tilde{q}}
\hyphenation{pa-ra-me-ter}
\topmatter
\title
Weak convergence of orthogonal polynomials
\endtitle
\author
Walter Van Assche
\endauthor
\affil
Katholieke Universiteit Leuven
\endaffil
\address
Department of Mathematics, Katholieke Universiteit Leuven,
Celestijnenlaan 200\,B, B-3001 Heverlee (BELGIUM)
\endaddress
\email
fgaee03\@cc1.kuleuven.ac.be
\endemail
\thanks
The author is a Senior Research Associate of the Belgian National Fund
for Scientific Research
\endthanks
\rightheadtext{Orthogonal Polynomials}

\abstract
The weak convergence of orthogonal polynomials is given under
conditions on the asymptotic behaviour of the coefficients in the three-term
recurrence relation. The results generalize known results and are applied
to several systems of orthogonal polynomials, including orthogonal polynomials
on a finite set of points.
 \endabstract

\endtopmatter

\document
\head 1. Introduction \endhead
Let $p_n(x)$ be a system of orthonormal polynomials on the real line,
with orthogonality measure $\mu$, i.e., $\mu$ is a probability measure
for which all the moments exist and
$$  \int p_n(x)p_m(x) \, d\mu(x) = \delta_{m,n}, \qquad m,n \geq 0 . $$
When the support of $\mu$ consists of a finite number of points
$x_1,x_2,\ldots,x_N$, then we only consider the polynomials up to degree
$N$ and $p_N(x)$ has its zeros at the support $\{x_1,x_2,\ldots,x_N\}$.
It is well known that orthonormal polynomials satisfy a three-term
recurrence relation
$$  xp_n(x) = a_{n+1}p_{n+1}(x) + b_np_n(x) + a_n p_{n-1}(x), \qquad
   n \geq 0, \tag 1.1 $$
with initial values $p_0(x) = 1$ and $p_{-1}(x) = 0$. Here
$$   a_{n+1} = \int x p_{n}(x)p_{n+1}(x) \, d\mu(x) ,\quad
     b_n = \int x p_n^2(x)\, d\mu(x) \in {\Bbb R}, \qquad n \geq 0 . $$
Usually the orthonormal polynomials are chosen in such a way that the
leading coefficient $\gamma_n = (a_1a_2\cdots a_n)^{-1}$ is positive, and then
$a_n =\gamma_n/\gamma_{n-1} > 0$ for every $n \geq 1$.

We are interested in the weak asymptotic behaviour of the orthonormal
polynomials $p_n(x)$. This means that we want to investigate the behaviour
for $n \to \infty$ of integrals of the form
$$   \int f(x) p_n^2(x) \, d\mu(x), \qquad f \in C_b , $$
where $C_b$ is the linear space of bounded and continuous functions
on ${\Bbb R}$.

We will consider a one-parameter orthogonality measure $\mu_k$ $(k \in {\Bbb
N})$, the parameter being discrete. This implies that the recurrence
coefficients and the orthogonal polynomials all depend on this discrete
parameter and we write
$$   xp_n(x;\mu_k) = a_{n+1,k}p_{n+1}(x;\mu_k) + b_{n,k}p_n(x;\mu_k)
   + a_{n,k} p_{n-1}(x;\mu_k), \qquad n\geq 0 . \tag 1.2 $$
Our main result will be the limit as $n \to \infty$ of integrals of the form
$$   \int f(x) p_{n+k}(x;\mu_n) p_{n+l}(x;\mu_n) \, d\mu_n(x), \qquad f \in C_b
. $$
The special case $k=l=0$ then gives the desired weak convergence, but the
general case with $k,l \in {\Bbb Z}$ also has useful applications: these
integrals are related to the transition probabilities of birth and death
processes and random walks.

Such one-parameter families do occur frequently in various applications
and limiting procedures. For example the rescaling of orthogonal polynomials
$$  p_n(c_k x;\mu) = p_n(x;\mu_k) $$
gives the one-parameter family of measures $\mu_k$ with distribution functions
satisfying
$$  \mu_k(t) = \mu(c_kt), \qquad t \in {\Bbb R}.  $$
Other examples include orthogonal polynomials in which some of the parameters
are allowed to tend to infinity together with the degree,
e.g., if $P_n^{(\alpha,\beta)}(x)$ is the Jacobi polynomial of degree $n$,
with weight function
$$   w(x) = (1-x)^\alpha (1+x)^\beta, \qquad -1 < x < 1, $$
then $P_n^{(k\alpha+\gamma,k\beta+\delta)}(x)$ is an orthogonal polynomial
of degree $n$ with weight function
$$  w_k(x) = w^k(x) (1-x)^\gamma (1+x)^\delta, \qquad -1 < x < 1. $$

The main result will be in terms of a doubly infinite Jacobi matrix
$$   \J = \pmatrix
   \ddots   & \ddots   &          &          &       &        &       & \\
   \ddots   & b_{-2}^0 & a_{-1}^0 &          &       &        &       & \\
            & a_{-1}^0 & b_{-1}^0 & a_{0}^0  &       &        &       & \\
            &          & a_0^0    & b_0^0    & a_1^0 &        &       & \\
            &          &          & a_1^0    & b_1^0 & a_2^0  &       & \\
            &          &          &          & a_2^0 & b_2^0  & a_3^0 & \\
            &        &        &        &     & \ddots & \ddots& \ddots
           \endpmatrix ,  $$
and is the following

\proclaim{Theorem}
Suppose that the recurrence coefficients in \thetag{1.2} satisfy
$$  \lim_{n \to \infty} a_{n+k,n} = a_k^0 > 0, \qquad
   \lim_{n \to \infty} b_{n+k,n} = b_k^0 \in {\Bbb R} \tag 1.3 $$
for every $k \in {\Bbb Z}$. Then
$$  \multline
   \lim_{n \to \infty}
\int f(x) p_{n+k}(x;\mu_n) p_{n+l}(x;\mu_n) \, d\mu_n(x)  \\
  = \int f(x)    \pmatrix A_k(x) & B_k(x) \endpmatrix d\mu(x)
                  \pmatrix A_l(x) \\ B_l(x) \endpmatrix,
  \endmultline                  $$
for every polynomial $f$.
Here $\mu$ is the spectral matrix of measures for the doubly
infinite Jacobi matrix $\J$ with entries $a_n^0, b_n^0$ $(n \in {\Bbb Z})$
and
$$  \pmatrix A_n(x) \\ B_n(x) \endpmatrix = \cases
    \pmatrix - \frac{a_0^0}{a_1^0} p_{n-1}^{(1)}(x) \\
                   p_n(x)   \endpmatrix &\text{for $n \geq 0$,} \\
  \pmatrix  q_{|n|-1}(x) \\ -\frac{a_0^0}{a_{-1}^0} q_{|n|-2}^{(1)}(x)
                  \endpmatrix&\text{for $n < 0$}, \endcases $$
with $p_n(x)$ and $q_n(x)$ the orthonormal polynomials with recurrence
coefficients respectively $a_{n+1}^0, b_n^0$ $(n=0,1,2,\ldots)$
and $a_{-n-1}^0, b_{-n-1}^0$, $(n=0,1,2,\ldots)$.
If the doubly infinite Jacobi matrix $\J$ has a unique self-adjoint extension
to its maximal domain $\{ \psi \in \ell_2({\Bbb Z},{\Bbb C}) :
\J \psi \in \ell_2({\Bbb Z},{\Bbb C}) \}$,
then the result holds also for  every bounded continuous
function $f$ on the real line.
\endproclaim

The proof of this theorem relies on various elements from the spectral
theory of Jacobi matrices. We will briefly review some results
from this spectral theory in Section 2 and then prove the theorem
in Section 3. In Section 4 we will apply  the theorem to some specific
systems of orthogonal polynomials.

\head 2. Spectral theory for Jacobi operators \endhead
If we put the coefficients $a_n>0$ $(n=1,2,\ldots)$ and $b_n\in {\Bbb R}$
$(n=0,1,2,\ldots)$ of the recurrence relation \thetag{1.1}
in a tridiagonal matrix
$$   J = \pmatrix
     b_0 & a_1 &     &        &  & \\
     a_1 & b_1 & a_2 &        &  & \\
         & a_2 & b_2 & a_3    &  & \\
         &     & a_3 & b_3    &  a_4 & \\
         &     &     &\ddots  & \ddots & \ddots
        \endpmatrix    \tag 2.1 $$
then $J$ is the Jacobi matrix associated with these orthogonal polynomials.
If $\mu$ is supported on $N$ points, then $a_n=0=b_n$ for $n\geq N$ and
$J$ is a $N\times N$ matrix with eigenvalues at the support
$\{x_1,x_2,\ldots,x_N\}$. Note that
$$   J = \left( \int xp_i(x)p_j(x) \,d\mu(x) \right)_{i,j=0,1,\ldots},
     \tag 2.2$$
and by induction we find
$$  J^k =  \left( \int x^kp_i(x)p_j(x) \,d\mu(x) \right)_{i,j=0,1,\ldots}.
     \tag 2.3 $$
In fact, $J: \ell_2({\Bbb N},{\Bbb C}) \to \ell_2({\Bbb N},{\Bbb C})$ acts as a
linear operator on the Hilbert
space $\ell_2({\Bbb N},{\Bbb C}) = \{ \psi : \psi_i \in {\Bbb C} \text{ and }
\sum_{i=0}^\infty |\psi_i|^2 < \infty\}$. The operator is symmetric
on the initial domain consisting of finite linear combinations
of the basic vectors
$\{ e_n^+ = (\underbrace{0,0,\ldots,0}_{\text{$n$ zeros}},1,0,\ldots):
n \geq 0 \}$
and by imposing appropriate conditions on the recurrence
coefficients $a_{n+1}, b_n$ (e.g., boundedness) this operator can be extended in
a unique way to a self-adjoint operator on the maximal domain
$\{ \psi \in \ell_2({\Bbb N},{\Bbb C}) : J\psi \in \ell_2({\Bbb N},{\Bbb C})
\}$. This operator has a cyclic
vector, i.e., if we take $e^+_0 = (1,0,0,0,\ldots)$, then the linear span of
$\{J^k e^+_0: k=0,1,\ldots\}$ is dense in $\ell_2({\Bbb N},{\Bbb C})$.
The spectral theorem (see, e.g., Akhiezer and Glazman \cite{2} or Stone
\cite{19}) then implies the existence of a measure $\mu$ and
a linear mapping
$\Lambda: \ell_2({\Bbb N},{\Bbb C}) \to L_2(\mu)$ with $\Lambda e^+_0 =1$
and $\Lambda J \psi = M \Lambda \psi$ for every
$\psi \in \ell_2({\Bbb N},{\Bbb C})$, where
$M$ is the multiplication operator $$   (M f)(t) = t f(t). $$
The mapping $\Lambda$ is unitary, meaning $\langle \Lambda \psi, \Lambda \phi
\rangle = \langle \psi,\phi \rangle$. The mapping $\Lambda$ thus maps $e^+_0$
to the constant function $t \mapsto 1$, $J e^+_0$ to the identity
$t \mapsto t$, and in general
$\Lambda$ maps $J^n e^+_0$ to the monomial $t \mapsto t^n$. Hence the
fact that $e^+_0$
is a cyclic vector is equivalent with the density of polynomials in $L_2(\mu)$.
By induction and by
using the recurrence relation \thetag{1.1} we see that $\Lambda$
maps the basic vector $e^+_n$ to the polynomial $p_n$. The unitarity thus
implies that
$$  \int p_n(t) p_m(t) \,d\mu(t) = \langle e^+_n,e^+_m \rangle = \delta_{m,n}, $$
which shows that the spectral measure $\mu$ for the operator $J$ is
the orthogonality measure for the orthogonal polynomials $p_n(x)$
$(n=0,1,2,\ldots)$. These elements from spectral theory, applied to the
semi-infinite Jacobi matrix $J$, are well-known (see e.g., Akhiezer \cite{1},
Dombrowski \cite{6}, Sarason \cite{18}, Stone \cite{19}).

We will also need to use doubly infinite Jacobi matrices of the form
$$   \J = \pmatrix
   \ddots   & \ddots   &      &        &     &        &       & \\
   \ddots   & b_{-2} & a_{-1} &        &     &        &       & \\
            & a_{-1} & b_{-1} & a_{0}  &     &        &       & \\
            &        & a_0    & b_0    & a_1 &        &       & \\
            &        &        & a_1    & b_1 & a_2    &       & \\
            &        &        &        & a_2 & b_2    & a_3   & \\
            &        &        &        &     & \ddots & \ddots& \ddots
           \endpmatrix , \tag 2.4 $$
with $a_k > 0, b_k \in {\Bbb R}$ for $k \in {\Bbb Z}$.
Again these are operators, now acting on the Hilbert space
$\ell_2({\Bbb Z},{\Bbb C}) = \{ \psi : \psi_i \in {\Bbb C} \text{ and }
\sum_{i=-\infty}^\infty |\psi_i|^2 < \infty\}$.
The spectral theory of such matrices is less known, but has also been
developped (see, e.g., Berezanski\u\i\ \cite{4}, Nikishin \cite{17},
Masson and Repka \cite{13}). We will
briefly
recall some of the elements, which we have taken from Nikishin \cite{17}
(see also Berezanski\u\i\ \cite{4, Chapter VII, \S 3}).
We will first consider bounded matrices $\J$, i.e., we assume that
$a_n$ and $b_n$ $(n \in {\Bbb Z})$ are bounded sequences. It turns out
that it is convenient to introduce the $2\times 2$ matrices
$$   B_0 = \pmatrix b_{-1} & a_0 \\
                    a_0    & b_0   \endpmatrix , \quad
B_n = \pmatrix  b_{-n-1} & 0 \\
                 0       & b_n  \endpmatrix  \qquad n=1,2,\ldots $$

$$   A_n = \pmatrix a_{-n} & 0 \\
                    0 & a_n  \endpmatrix , \qquad n=1,2,\ldots. $$
We can now study the semi-infinite Jacobi block matrix
$$   {\Bbb J} = \pmatrix
     B_0 & A_1 &     &        &  & \\
     A_1 & B_1 & A_2 &        &  & \\
         & A_2 & B_2 & A_3    &  & \\
         &     & A_3 & B_3    &  A_4 & \\
         &     &     &\ddots  & \ddots & \ddots
        \endpmatrix    \tag 2.5 $$
which contains $2\times 2$ matrices. As an operator it acts on
$\ell_2({\Bbb N},{\Bbb C}^2)$ in the sense that
$\J \psi$ for $\psi \in \ell_2({\Bbb Z},{\Bbb C})$
corresponds to ${\Bbb J} \Psi$ with $\Psi \in \ell_2({\Bbb N},{\Bbb C}^2)$
given by
$$ \Psi_n = \pmatrix \psi_{-n-1} \\ \psi_{n} \endpmatrix,
    \qquad n=0,1,2,\ldots. $$
In this way we have transformed the study of the doubly-infinite Jacobi matrix
to the study of a semi-infinite Jacobi block matrix. Such semi-infinite
block matrices are closely connected to orthogonal matrix polynomials,
in the same way as ordinary Jacobi matrices are connected to scalar
orthogonal polynomials (see, e.g., Aptekarev and Nikishin \cite{3}).

Consider the standard basis $e_n$  $(n \in {\Bbb Z})$
in $\ell_2({\Bbb Z},{\Bbb C})$, i.e.,
$$  (e_n)_i = \delta_{i,n}, \qquad   i,n \in {\Bbb Z}, $$
then the linear span of $\{\J^k e_{-1}, \J^l e_0, k,l \in {\Bbb N}\}$ is dense
in $\ell_2({\Bbb Z},{\Bbb C})$. By the spectral theorem
there exists a matrix-measure
$$   \mu = \pmatrix \mu_{1,1}  & \mu_{1,2} \\
                    \mu_{1,2} & \mu_{2,2} \endpmatrix  $$
and a unitary linear mapping $\Lambda: \ell_2({\Bbb Z},{\Bbb C}) \to L_2(\mu)$
with
$$  \Lambda e_{-1} = \pmatrix 1 \\ 0 \endpmatrix, \quad
  \Lambda e_0 = \pmatrix 0 \\ 1 \endpmatrix, $$
such that
$$   \Lambda \J \psi = M \Lambda \psi, $$
where now $M$ is the multiplication operator in the space
$L_2(\mu)$ of vector valued functions. The inner product in the space
$L_2(\mu)$  is given by
$$   \langle \pmatrix f_1 \\ f_2 \endpmatrix ,
             \pmatrix g_1 \\ g_2 \endpmatrix \rangle =
    \int \pmatrix \overline{f_1(x)} & \overline{f_2(x)} \endpmatrix
        \pmatrix d\mu_{1,1}(x) & d\mu_{1,2}(x) \\
                 d\mu_{1,2}(x) & d\mu_{2,2}(x) \endpmatrix
    \pmatrix g_1(x) \\ g_2(x) \endpmatrix . $$
Therefore the mapping $\Lambda$ maps $\J^n e_{-1}$ to the vector function
$$      t \mapsto \pmatrix t^n \\ 0 \endpmatrix , $$
whereas $\J^n e_0$ is mapped to
$$      t \mapsto \pmatrix 0 \\ t^n \endpmatrix . $$

Let $J^+$ be the semi-infinite Jacobi matrix
$$   J^+ = \pmatrix
             b_0    & a_1 &        &       &  \\
             a_1    & b_1 & a_2    &       &  \\
                    & a_2 & b_2    & a_3   &  \\
                    &     & \ddots & \ddots& \ddots
           \endpmatrix ,$$
and similarly $J^-$ be the semi-infinite Jacobi matrix
$$   J^- = \pmatrix
              b_{-1} & a_{-1} &        &        & \\
              a_{-1} & b_{-2} & a_{-2} &        & \\
                     & a_{-2} & b_{-3} & a_{-3} & \\
                     &        & \ddots & \ddots & \ddots
           \endpmatrix.  $$
We denote by $p_n(x)$ the orthonormal polynomials corresponding to the
Jacobi-matrix $J^+$ satisfying the recurrence relation
$$  xp_n(x) = a_{n+1}p_{n+1}(x) + b_np_n(x) + a_n p_{n-1}(x),
   \tag 2.6 $$
with initial values $p_{-1}(x)=0, p_0(x)=1$, and by $q_n(x)$ the
orthonormal polynomials
for $J^-$ satisfying
$$  xq_n(x) = a_{-n-1}q_{n+1}(x) + b_{-n-1}q_n(x) + a_{-n} q_{n-1}(x),
   \tag 2.7 $$
with initial values $q_{-1}(x)=0, q_0(x)=1$.
We will show, by induction, that for $n \in {\Bbb N}$ the mapping
$\Lambda$ maps the basis vector $e_n$ to $\Lambda e_n$ which is the vector
function
$$ t \mapsto \pmatrix      - \frac{a_0}{a_1} p_{n-1}^{(1)}(t) \\
                   p_n(t)   \endpmatrix  \tag 2.8 $$
and the basisvector $e_{-n}$ to $\Lambda e_{-n}$ given by
$$ t \mapsto \pmatrix       q_{n-1}(t) \\
           -\frac{a_0}{a_{-1}} q_{n-2}^{(1)}(t)   \endpmatrix . \tag 2.9 $$
Here $p_n^{(1)}(x)$ and $q_n^{(1)}(x)$ are the associated polynomials, i.e.,
the orthonormal polynomials
corresponding to the Jacobi matrices $J^+$ and $J^-$, with the first row and
column deleted.
This is clear for $e_{-1}$ and $e_0$. Assume that this is true for $0 \leq n
\leq k$, then from
$$   \J e_k = a_{k+1} e_{k+1} + b_k e_{k} + a_k e_{k-1}  $$
it follows that
$$  \Lambda \J e_k = a_{k+1} \Lambda e_{k+1} + b_k \Lambda e_{k} + a_k
\Lambda e_{k-1}, $$
and since $\Lambda \J e_k = M \Lambda e_k$ this gives
$$     (\Lambda e_{k+1})(t) = \frac{1}{a_{k+1}}
         \pmatrix
         -\frac{a_0}{a_1} \left( (t-b_k)p_{k-1}^{(1)}(t) - a_kp_{k-2}^{(1)}(t)
         \right) \\
          (t-b_k) p_k(t) - a_k p_{k-1}(t)
         \endpmatrix
         = \pmatrix       - \frac{a_0}{a_1} p_k^{(1)}(t) \\
                   p_{k+1}(t)   \endpmatrix.  $$
Similarly for $e_{-n}$ the result is true for $n=0$ and $n=1$ and since
$$   \J e_{-k} = a_{-k+1} e_{-k+1} + b_{-k} e_{-k} + a_{-k} e_{-k-1} $$
and $\Lambda \J e_{-k} = M \Lambda e_{-k}$, this gives
$$  (\Lambda e_{-k-1})(t) = \frac{1}{a_{-k}}
         \pmatrix
          (t-b_{-k})q_{k-1}(t) - a_{-k+1} q_{k-2}(t) \\
         -\frac{a_0}{a_{-1}} \left( (t-b_{-k}) q_{k-2}^{(1)}(t) - a_{-k+1}
        q_{k-3}^{(1)}(t) \right) \endpmatrix
         = \pmatrix      q_k(t) \\
                   -\frac{a_0}{a_{-1}} q_{k-1}^{(1)}(t) \endpmatrix.  $$
From the unitarity we find for $m,n \in {\Bbb Z}$
$$  \langle \Lambda e_n, \Lambda e_m \rangle =  \delta_{m,n}, $$
Hence the matrix polynomials
$$  P_n(t) = \pmatrix q_n(t) & -\frac{a_0}{a_{-1}} q_{n-1}^{(1)}(t)   \\
           -\frac{a_0}{a_1} p_{n-1}^{(1)}(t) & p_{n}(t) \endpmatrix $$
satisfy
$$    \int P_n(t) d\mu(t) P_m(t)^* =  \pmatrix 1 & 0 \\ 0 & 1 \endpmatrix
   \delta_{m,n} . $$
Therefore these matrix polynomials are orthonormal with respect to the
matrix-measure $\mu$.
Note that these matrix polynomials satisfy
$$  \align
t P_n(t) &= A_{n+1} P_{n+1}(t) + B_n P_n(t) + A_n P_{n-1}(t), \\
  P_0(t) &= \pmatrix  1 & 0 \\ 0 & 1 \endpmatrix ,
  \endalign $$
so that they are the orthonormal polynomials corresponding to the
block Jacobi matrix ${\Bbb J}$ given in \thetag{2.5}.

The polynomials $p_n(x)$ $(n=0,1,2,\ldots)$ are orthonormal with respect
to some probability measure $\mu^+$, and the polynomials $q_n(x)$
$(n=0,1,2,\ldots)$ are orthonormal with respect to some probability measure
$\mu^-$. In order to find a relation between the matrix-measure $\mu$
and the measures $\mu^+$ and $\mu^-$, we observe that the unitarity
of $\Lambda$ implies
$$    \langle (z-\J)^{-1} e_n , e_m \rangle =
   \int \frac{1}{z-t} (\Lambda e_n)^* d\mu(t) \Lambda e_m, $$
where $(z-\J)^{-1}$ is the resolvent of $\J$, which is well-defined
for every $z$ outside the spectrum of $\J$. Since $\J$ is symmetric
and bounded, it follows that $\J$ is self-adjoint so that its spectrum
is a subset of the real line. The Stieltjes transform of the matrix-measure
$\mu$ is determined by
$$ \langle (z-\J)^{-1} e_{-1} , e_{-1} \rangle =
   \int \frac{1}{z-t} d\mu_{1,1}(t), \quad
    \langle (z-\J)^{-1} e_{-1} , e_0 \rangle =
   \int \frac{1}{z-t} d\mu_{1,2}(t), $$
$$   \langle (z-\J)^{-1} e_0 , e_0 \rangle =
   \int \frac{1}{z-t} d\mu_{2,2}(t).
  $$
On the other hand, if we write
$$   (z-\J)^{-1} e_0 = r = (\ldots, r_{-2}, r_{-1}, r_0, r_1, r_2, \ldots), $$
then   $(z-\J) r = e_0$, which gives the infinite system of equations
$$ \align
 z r_k &= a_k r_{k-1} + b_k r_k + a_{k+1} r_{k+1}, \qquad k \geq 1, \\
 z r_{-k} &= a_{-k} r_{-k-1} + b_{-k} r_{-k} + a_{-k+1} r_{-k+1}, \qquad k \geq
        1, \tag 2.10   \\
  zr_0 - 1 &= a_0 r_{-1} + b_0 r_0 + a_1 r_1 .
  \endalign  $$
These recurrence relations for the unknowns $r_k$ are precisely the
three-term recurrence relations associated with the Jacobi-matrices
$J^+$ and $J^-$. The general solution of the equations
can be expressed as
a linear combination
of the orthogonal polynomials $p_n(z)$ (respectively $q_n(z)$) and the
functions of the second kind $\p_n(z)$ (respectively $\q_n(z)$)  given by
$$  \p_n(z) = \int \frac{p_n(x)}{z-x} \, d\mu^+(x), \quad
  \q_n(z) = \int \frac{q_n(x)}{z-x} \, d\mu^-(x), $$
with $a_0\p_{-1}(z) = 1 = a_0 \q_{-1}(z)$.
These functions of the second kind have the property that they are a minimal
solution of the recurrence relation, and for $z \in {\Bbb C} \setminus {\Bbb R}$
they satisfy $\lim_{n \to \infty} \p_n(z) = \lim_{n \to \infty} \q_n(z) = 0$.
The fact that $r \in \ell_2({\Bbb Z},{\Bbb C})$ thus implies that
$r_n$ and $r_{-n}$ are (up to a constant factor) given by the functions
of the second kind $\p_n(z)$ and $\q_{n-1}(z)$ respectively. The constant
multiple is determined by setting $n=0$, giving
$$  r_n = r_0 \frac{\p_{n}(z)}{\p_0(z)} , \quad r_{-n} = a_0r_0 \q_{n-1}(z),
\qquad n \geq 0. $$
In particular $r_1 = r_0 \p_1(z)/\p_0(z)$ and $r_{-1} = a_0 r_0 \q_0(z)$.
Inserting in \thetag{2.10} gives
$$    r_0 = \frac{1}{z - a_0^2 \q_0(z) - b_0 - a_1 \p_1(z)/\p_0(z)} , $$
which by using $a_1\p_1(z) = (z-b_0)\p_0(z) -1$ gives
$$    r_0 = \frac{\p_0(z)}{1-a_0^2 \p_0(z)\q_0(z)} . $$
In a similar way we may investigate
$$   (z-\J)^{-1} e_{-1} = s = (\ldots, s_{-2}, s_{-1}, s_0, s_1, s_2, \ldots),
$$
which gives the infinite system of linear equations
$$ \align
 z s_k &= a_k s_{k-1} + b_k s_k + a_{k+1} s_{k+1}, \qquad k \geq 0, \\
 z s_{-k} &= a_{-k} s_{-k-1} + b_{-k} s_{-k} + a_{-k+1} s_{-k+1}, \qquad k \geq
       2,  \tag 2.11  \\
  zs_{-1} - 1 &= a_{-1} s_{-2} + b_{-1} s_{-1} + a_0 s_0 .
  \endalign  $$
Now we find
$$  s_n = a_0 s_{-1} \p_n(z), \quad s_{-n} = s_{-1}
         \frac{\q_{n-1}(z)}{\q_0(z)}, $$
and inserting this in \thetag{2.11} gives
$$   s_{-1} = \frac{1}{z-a_{-1} \q_{1}(z)/\q_0(z) - b_{-1} - a_0^2 \p_0(z)}, $$
which by using $a_{-1}\q_1(z) = (z-b_{-1})\q_0(z) - 1$ becomes
$$  s_{-1} = \frac{\q_0(z)}{1-a_0^2 \p_0(z)\q_0(z)} . $$
The Stieltjes transform of the matrix of measures $\mu$ is thus given by
$$ \int \frac{1}{z-t} d\mu_{1,1}(t) = \frac{\q_0(z)}{1-a_0^2 \p_0(z)\q_0(z)},
  \quad \int \frac{1}{z-t} d\mu_{2,2}(t) =
\frac{\p_0(z)}{1-a_0^2\p_0(z)\q_0(z)}, $$
$$ \int \frac{1}{z-t} d\mu_{1,2}(t) = \frac{a_0
\p_0(z)\q_0(z)}{1-a_0^2\p_0(z)\q_0(z)}, $$
where $\p_0(z)$ and $\q_0(z)$ are respectively the Stieltjes transforms of
the measures $\mu^+$ and $\mu^-$.

So far we have supposed that $\J$ is a bounded operator, but all the above
results
remain true whenever $\J$ has a unique self-adjoint extension
to its maximal domain $\{ \psi \in \ell_2({\Bbb Z},{\Bbb C}) :
\J \psi \in \ell_2({\Bbb Z},{\Bbb C}) \}$,  which
will be true whenever both $J^+$ and $J^-$ are Jacobi matrices
corresponding to a determinate moment problem.

\head 3. Proof of the theorem \endhead

The first important observation is that for $f(x)=x^m$ one has
$$  \int x^m p_{n+k}(x;\mu_n) p_{n+l}(x;\mu_n) \, d\mu_n(x) =
    \langle J_n^m e^+_{n+k},  e^+_{n+l} \rangle , $$
where $J_n$ is the semi-infinite Jacobi matrix with the recurrence
coefficients $a_{k,n}, b_{k,n}$, $(k = 0,1,2,\ldots)$.
Consider the semi-infinite operator $J_n$ as a doubly infinite
Jacobi matrix by taking $(J_n)_{i,j} = 0$ whenever $i<0$ or $j<0$.
If $S$ is the shift operator on $\ell_2({\Bbb Z},{\Bbb C})$ acting as
$S e_k = e_{k+1}$, then
$$  \int x^m p_{n+k}(x;\mu_n) p_{n+l}(x;\mu_n) \, d\mu_n(x) =
         \langle (S^*)^n J_n^m S^n e_k,e_l \rangle . $$
One easily verifies that this expression is given by
$$ \multline
\langle (S^*)^n J_n^m S^n e_k,e_l \rangle
   = \sum\Sb i_1,i_2,\ldots,i_{m-1} \in \{-1,0,1\}\\i_1+i_2+\cdots+i_{m-1}+k-l
     \,\in \{-1,0,1\}\endSb
     (J_n)_{n+k,n+k+i_1} (J_n)_{n+k+i_1,n+k+i_1+i_2} \\
     \cdots (J_n)_{n+k+i_1+i_2+\cdots+i_{m-1},n+l}.
   \endmultline  $$
This is a finite sum, containing the matrix entries $(J_n)_{n+r,n+s}$
where $r$ and $s$ remain bounded.
The hypothesis \thetag{1.3} implies that $(S^*)^n J_n S^n$ converges
entrywise to $\J$,
where $\J$ is the doubly infinite Jacobi matrix containing the coefficients
$a_n^0,b_n^0$, $(n \in {\Bbb Z})$.
The hypothesis \thetag{1.3} thus implies that
$$ \multline
\lim_{n \to \infty} \int x^m p_{n+k}(x;\mu_n) p_{n+l}(x;\mu_n) \, d\mu_n(x)\\
   = \sum\Sb i_1,i_2,\ldots,i_{m-1} \in \{-1,0,1\}\\i_1+i_2+\cdots+i_{m-1}+k-l
     \,\in \{-1,0,1\}\endSb
      \J_{k,k+i_1} \J_{k+i_1,k+i_1+i_2} \cdots
     \J_{k+i_1+i_2+\cdots+i_{m-1},l},
  \endmultline   $$
and the latter expression is equal to
$$   \langle \J^m e_k , e_l \rangle . $$
Therefore we find that $(S^*)^n J_n^m S^n$ converges entrywise to $\J^m$.
By the unitarity of the mapping $\Lambda: \ell_2({\Bbb Z},{\Bbb C}) \to
L_2(\mu)$, transforming the action of $\J$ on $\ell_2({\Bbb Z},{\Bbb C})$ to the
action of the multiplication operator $M$ on $L_2(\mu)$,
we have
$$ \langle \J^m e_k , e_l \rangle  =
\int t^m (\Lambda e_k)^* d\mu(t) \Lambda e_l, $$
and thus the theorem for $f(x)=x^m$ (and hence for every polynomial $f$)
follows from \thetag{2.8} and \thetag{2.9}.

In order to proof the theorem for every $f \in C_b$, we consider
the linear operator $H_n = (S^*)^n J_n S^n$.
We would like to prove that $f(H_n)=(S^*)^n f(J_n) S^n$
converges weakly to $f(\J)$, since then
$$ \lim_{n \to \infty} \langle (S^*)^n f(J_n) S^n e_k,e_l \rangle
   = \langle f(\J) e_k, e_l \rangle, $$
and this is precisely the weak convergence stated in our theorem.
Observe that $H_n e_k = a_{k+n,n} e_{k-1} + b_{k+n,n} e_k + a_{k+n+1,n}
e_{k+1}$ whenever $k+n > 0$, hence
the condition \thetag{1.3} implies the convergence
of $H_n e_k$ to $\J e_k = a_k^0 e_{k-1} + b_k^0 e_k + a_{k+1}^0 e_{k+1}$
for every $k \in {\Bbb Z}$. All finite linear
combinations of the basis elements $e_k$ form a core  of $\J$ because
we assume that there is only one self-adjoint extension of $\J$ to
its maximal domain. Each $\psi$ in this core belongs to the domain of every
$H_n$ and
 $H_n \psi \to \J \psi$, hence $H_n$ converges to $\J$ in the sense of strong
 resolvent convergence, i.e.,
 for all $z \in {\Bbb C} \setminus {\Bbb R}$ we have strong
 convergence of $(z-H_n)^{-1}$ to $(z-\J)^{-1}$
 \cite{22, Theorem 9.16 (i) on p.~283}. But then the strong resolvent
 convergence of $H_n$ to $\J$ implies that $f(H_n)$ converges strongly
to $f(\J)$ for every bounded and continuous function on ${\Bbb R}$
\cite{22, Theorem 9.17 on p.~ 284}, and
strong convergence in turn implies weak convergence and in particular
the convergence of $\langle (S^*)^n f(J_n) S^n e_k,e_l \rangle$ to
$\langle f(\J) e_k,e_l \rangle$, which is the desired asymptotic behaviour.

\head 4. Examples \endhead

\subhead The class $M(a,b)$ \endsubhead
The class $M(a,b)$ consists of all orthogonal polynomials $p_n(x;\mu)$ (or all
probability measures $\mu$)  with recurrence coefficients that satisfy
$$    \lim_{n \to \infty} a_n = a/2, \quad \lim_{n \to \infty} b_n = b. $$
We can apply the theorem with the family of measures $\mu_k \equiv \mu$,
i.e., with all the orthogonality measures the same. If $a>0$ then we can,
without loss of generality, only consider $M(1,0)$. The doubly infinite
Jacobi matrix $\J$ then consists of $0$ on the diagonal and $1/2$ on the
subdiagonals. The semi-infinite Jacobi matrices $J^+$ and $J^-$ are
the same and the corresponding orthogonal polynomials are the Chebyshev
polynomials of the second kind, which are orthogonal with
respect to the measure $(2/\pi) \sqrt{1-x^2}\, dx$ on the interval
$[-1,1]$. The Stieltjes transform of this measure is
$$  \p_0(z) = \q_0(z) = 2 [ z - \sqrt{z^2-1} ]. $$
Therefore the Stieltjes transform of the spectral matrix of measures for the
Jacobi matrix $\J$ is given by
$$   \int \frac{1}{z-x} \, d\mu_{1,1}(x) = \frac{1}{\sqrt{z^2-1}} =
      \int \frac{1}{z-x} \, d\mu_{2,2}(x), $$
$$ \int \frac{1}{z-x} \, d\mu_{1,2}(x) = \frac{z-\sqrt{z^2-1}}{\sqrt{z^2-1}},
$$
from which one easily finds that the spectrum of $\J$ is $[-1,1]$ and
$$  d\mu_{1,1}(x) = d\mu_{2,2}(x) = \frac{1}{\pi} \frac{dx}{\sqrt{1-x^2}} ,
\quad d\mu_{1,2}(x) = \frac{1}{\pi} \frac{x\, dx}{\sqrt{1-x^2}}. $$
From our theorem we thus find
$$  \lim_{n \to \infty} \int f(x) p_n^2(x;\mu)\, d\mu(x)
   = \int f(x) \, d\mu_{2,2}(x) = \frac{1}{\pi}\int_{-1}^1
\frac{f(x)}{\sqrt{1-x^2}} \, dx, $$
and in general
$$ \align
\lim_{n \to \infty} \int f(x) p_n(x;\mu)p_{n+k}(x;\mu)\, d\mu(x)
  &= \frac{1}{\pi} \int_{-1}^1 f(x) \frac{U_{k}(x) - x
U_{k-1}(x)}{\sqrt{1-x^2}} \, dx  \\
  &= \frac{1}{\pi} \int_{-1}^1 f(x) \frac{T_k(x)}{\sqrt{1-x^2}} \, dx.
  \endalign $$
Here we have used the identity
$$   U_{k}(x) - x U_{k-1}(x) = T_k(x), $$
where $T_k(x)$ is the Chebyshev polynomial of the first kind.
This result is well-known and can already be found in \cite{15, Theorem 13
on p.\ 45}. See also \cite{20, Theorem 2 on p.\ 438}.

\subhead Unbounded recurrence coefficients \endsubhead
Suppose that we have a sequence of orthogonal polynomials $p_n(x;\mu)$
satisfying the three-term recurrence relation \thetag{1.1}. If we rescale the
variable by a positive and increasing sequence $c_k$ and consider the
one-parameter family of polynomials $p_n(c_kx;\mu)$, then these polynomials
satisfy a recurrence relation of the form \thetag{1.2} with
$$    a_{n,k} = \frac{a_n}{c_k}, \quad b_{n,k} = \frac{b_n}{c_k}. $$
Suppose the sequence $c_k$ is such that
$$  \lim_{n \to \infty} \frac{a_n}{c_n} = a/2 > 0, \quad
   \lim_{n \to \infty} \frac{b_n}{c_n} = b,  $$
and
$$   \lim_{n \to \infty} \frac{c_{n+1}}{c_n} = 1 , $$
then it is easy to show that
$$   \lim_{n \to \infty} a_{n+k,n} = a/2,  \quad
     \lim_{n \to \infty} b_{n+k,n} = b. $$
Therefore the conditions of our theorem are valid and the theorem is true,
with $\J$ the doubly infinite Jacobi matrix with constant entries $b$ on the
diagonal and $a/2$ on the subdiagonal. The orthogonal polynomials
for the Jacobi matrices $J^+$ and $J^-$ are $U_n(\frac{x-b}{a})$
and the matrix of measures is supported on $[b-a,b+a]$ and given by
$$  d\mu_{1,1}(x) = d\mu_{2,2}(x) = \frac{1}{\pi} \frac{dx}{\sqrt{a^2-(x-b)^2}}
, \quad d\mu_{1,2}(x) = \frac{1}{a\pi} \frac{(x-b)\, dx}{\sqrt{a^2-(x-b)^2}}. $$
We thus have
$$  \lim_{n \to \infty} \int f(x/c_n) p_n^2(x;\mu)\, d\mu(x)
   = \int f(x) \, d\mu_{2,2}(x) = \frac{1}{\pi}\int_{b-a}^{b+a}
\frac{f(x)}{\sqrt{a^2-(x-b)^2}} \, dx. $$
This result can already be found in \cite{10, Lemma 1 on p.\ 52} and
\cite{16, Lemma 3 on p.\ 1188}. In general we have
$$ \lim_{n \to \infty} \int f(x/c_n) p_n(x;\mu)p_{n+k}(x;\mu)\, d\mu(x)
  = \frac{1}{\pi} \int_{b-a}^{b+a} f(x) \frac{T_k((x-b)/a)}{\sqrt{a^2-(x-b)^2}}
\, dx.  $$

\subhead Wall polynomials \endsubhead
The orthonormal Wall polynomials $w_n(x;b,q)$, with $(0 < q < 1, 0 < b < 1)$ are
orthogonal on the geometric sequence $\{ q^n, n=1,2,3,\ldots\}$ and
have recurrence coefficients
$$   a_n(b,q) = q^n\sqrt{b(1-q^n)(1-bq^{n-1})}, \quad
     b_n(b,q) = q^n[b+q-(1+q)bq^n]. $$
If we consider the Wall polynomials $w_n(x;b;c^{1/k})$, where $0 < c < 1$,
then we have a one-parameter family of orthogonal polynomials with recurrence
coefficients
$$    a_{n,k} = a_n(b,c^{1/k}), \quad b_{n,k} = b_n(b,c^{1/k}), $$
and one easily verifies that
$$  \lim_{n \to \infty} a_{n+k,n} = c\sqrt{b(1-c)(1-bc)} = A/2, \quad
  \lim_{n \to \infty} b_{n+k,n} = (b+1-2bc)c = B. $$
Hence we can apply the theorem, where $\J$ again is a doubly infinite
Jacobi matrix with constant entries $B$ on the diagonal and $A/2$ on
the subdiagonals. Therefore again we have the asymptotic behaviour
in terms of Chebyshev polynomials (first and second kind) as in the previous
two examples. These three examples are all covered by Theorem 2 in
\cite{21, p.\ 307} which covers the special case of our theorem with
a doubly infinite Jacobi matrix with constant entries.
In particular this asymptotic behaviour for Wall polynomials was used
in \cite{21} to show that the product formulas for Legendre polynomials
are a limiting case of the product formulas for little $q$-Legendre
polynomials as $q \to 1$.

\subhead Jacobi polynomials $P_n^{(an+\alpha,bn+\beta)}(x)$
  \endsubhead
The Jacobi polynomials $P_n^{(\alpha,\beta)}(x)$ are orthogonal on
$[-1,1]$ with the weight function $(1-x)^\alpha(1+x)^\beta$. The orthonormal
Jacobi polynomials
$$  p_n^{(\alpha,\beta)}(x)  = \sqrt{\frac{2n+\alpha+\beta+1}{n+\alpha+\beta+1}
 \frac{n! (\alpha+\beta+2)_n}{(\alpha+1)_n(\beta+1)_n}} P_n^{(\alpha,\beta)}(x)
$$
have recurrence coefficients
$$ \align
a_n^2(\alpha,\beta) &= \frac{4n(n+\alpha)(n+\beta)(n+\alpha+\beta)}
  {(2n+\alpha+\beta-1)(2n+\alpha+\beta)^2(2n+\alpha+\beta+1)}  \\
 b_n(\alpha,\beta) &=
\frac{\beta^2-\alpha^2}{(2n+\alpha+\beta)(2n+\alpha+\beta+2)} .
  \endalign  $$
If we consider the Jacobi polynomials $p_n^{(ak+\alpha,bk+\beta)}(x)$
with $a>0, b >0$, then the
recurrence coefficients are $a_{n,k} = a_n(ak+\alpha,bk+\beta)$ and
$b_{n,k} = b_n(ak+\alpha,bk+\beta)$, and one easily finds
$$  \lim_{n \to \infty} a_{n+k,n} =
\frac{2\sqrt{(a+1)(b+1)(a+b+1)}}{(a+b+2)^2},
\quad  \lim_{n \to \infty} b_{n+k,n} =  \frac{b^2-a^2}{(a+b+2)^2}   . $$
Hence our theorem applies again with a constant doubly infinite Jacobi
matrix $\J$. This result has not been given in the literature, but
complements the known results concerning strong asymptotics and
zero behaviour given in  \cite{5} \cite{9} \cite{14}.

\subhead Laguerre polynomials $L_n^{an+\alpha}(nx)$ \endsubhead
The Laguerre polynomials $L_n^\alpha$ are orthogonal on $[0,\infty)$
with  weight function $x^\alpha e^{-x}$. The orthonormal Laguerre
polynomials
$$  p_n^\alpha(x) = (-1)^n \binom{n+\alpha}{n}^{-1/2} L_n^\alpha(x) $$
have recurrence coefficients
$$   a_n(\alpha) = \sqrt{n(n+\alpha)}, \quad b_n(\alpha) = 2n+\alpha+1. $$
If we consider the polynomials $p_n^{ak+\alpha}(kx)$, then
we have $a_{n,k} = a_n(ak+\alpha)/k$ and $b_{n,k} = b_n(ak+\alpha)/k$, hence
one easily finds
$$  \lim_{n \to \infty} a_{n+k,n} = a+1, \quad \lim_{n \to \infty}
    b_{n+k,n} = a+2, $$
so that once more our theorem can be applied with a constant Jacobi matrix
$\J$. This complements the asymptotic behaviour for such Laguerre
polynomials given in \cite{5} and \cite{8}.

\subhead Dual Hahn polynomials \endsubhead
The dual Hahn polynomials $R_n(x)=R_n(x;\alpha,\beta,N)$ are given by  the
recurrence relation
$$  -x R_k(x) = D_k R_{k-1}(x) - (D_k+B_k) R_k(x) + B_k R_{k+1}(x), $$
with initial condition $R_0(x)=1$ and $R_{-1}(x) = 0$ \cite{11}.
Here
$$   B_k = (N-1-k)(\alpha+1+k), \quad D_k = k(N+\beta-k) , $$
and the polynomials are orthogonal on the quadratic lattice
$\{x_k=k(k+\alpha+\beta+1): k=0,1,2,\ldots,N-1\}$ with weights
$$ \align
  \pi_k &= \pi_k(\alpha,\beta,N)     \\
   &=  \binom{N-1}{k} \frac{\Gamma(\beta+N)}{\Gamma(N+\alpha+\beta+k+1)}
      \frac{\Gamma(k+\alpha+1)\Gamma(k+\alpha+\beta+1)}{\Gamma(k+\beta+1)
      \Gamma(\alpha+1)} \, (2k+\alpha+\beta+1)
  \endalign    $$
at these points $x_k$,
so that these polynomials
are only defined up to degree $N$. The orthonormal polynomials are
$$  p_n(x;\alpha,\beta,N) = \frac{\binom{N+\alpha+\beta}{N-1}^{1/2}}
   {\binom{\alpha+n}{n}^{1/2} \binom{\beta+N-1-n}{N-1-n}^{1/2}} \
   R_n(x;\alpha,\beta,N), $$
with recurrence coefficients $a_n^2=D_nB_{n-1}$ and $b_n=D_n+B_n$.
These polynomials are useful in the description of a genetic  model
of Moran, as is worked out in \cite{11} \cite{12} and \cite{7}.
Consider the polynomials $p_n(k x;\alpha,\beta,k)$, then
we have the recurrence coefficients
$$ \align
a_{n,k}^2 &= \cases \frac{n(k-n)(k+\beta-n)(\alpha+n)}{k^2} &\text{for $n<k$}
     \\             0 & \text{for $n \geq k$}, \endcases \\
 b_{n,k} &= \cases \frac{(k-n-1)(\alpha+n+1) + n(k+\beta-n)}{k}
 &\text{for $n<k$} \\  0 & \text{for $n \geq k$}. \endcases
  \endalign $$
One easily finds
$$  \lim_{n \to \infty} a_{n+k,n}^2 = \cases 0 & \text{for $k \geq 0$}, \\
                                -k(\beta-k) & \text{for $k < 0$}, \endcases
   \quad
  \lim_{n \to \infty} b_{n+k,n} = \cases 0 & \text{for $k \geq 0$}, \\
                                -2k+\beta-1 & \text{for $k < 0$}. \endcases $$
The doubly infinite Jacobi matrix $\J$ corresponding to these asymptotic
formulas is therefore only a semi-infinite Jacobi matrix which coincides
with the Jacobi matrix $J^-$, for which the corresponding orthogonal
polynomials are the Laguerre polynomials $L_n^\beta(x)$. The spectral matrix
of measures of $\J$ thus reduces to $\mu_{2,2}\equiv 0 \equiv \mu_{1,2}$ and
$d\mu_{1,1}(x)=x^\beta e^{-x}/\Gamma(\beta+1)$
on $[0,\infty)$. Our theorem can easily be proved also for the case where
$\J$ reduces to a semi-infinite Jacobi matrix, and we thus have
$$  \lim_{n \to \infty}  \sum_{j=0}^{n-1}
    f(x_j/n) p_{n-1}^2(x_j;\alpha,\beta,n) \pi_{j,n}
   = \int_0^\infty f(x) \frac{x^\beta e^{-x}}{\Gamma(\beta+1)} \, dx, $$
where $x_j=j(j+\alpha+\beta+1)$ and $\pi_{j,n} = \pi_j(\alpha,\beta,n)$.
In general we have for $k,l \geq 1$
$$ \multline
\lim_{n \to \infty} \sum_{j=0}^{n-1}
  f(x_j/n) p_{n-k}(x_j;\alpha,\beta,n) p_{n-l}(x_j;\alpha,\beta,n) \pi_{j,n}
\\ = \frac{(-1)^{k+l}}{\sqrt{h_{k-1}h_{l-1}}}
 \int_0^\infty f(x) L_{k-1}^\beta(x) L_{l-1}^\beta(x) \frac{x^\beta
e^{-x}}{\Gamma(\beta+1)}\, dx, \endmultline $$
where $h_k = \binom{k+\beta}{k}$ is the norm of $L_k^\beta(x)$.

If we consider the polynomials  $p_n(k^{3/2}x+k^2/2;\alpha,k/2,k)$, then
the recurrence coefficients are zero for $n \geq k$ and for $n < k$ we have
$$ \align
a_{n,k}^2 &= \frac{n(k-n)(3k/2-n)(\alpha+n)}{k^3}, \\
 b_{n,k} &=  \frac{(k-n-1)(\alpha+n+1) + n(3k/2-n)-k^2/2}{k^{3/2}}.
\endalign $$
Now one easily finds
$$  \lim_{n \to \infty} a_{n+k,n}^2 = \cases 0 & \text{for $k \geq 0$}, \\
                                -k/2 & \text{for $k < 0$}, \endcases
   \quad
  \lim_{n \to \infty} b_{n+k,n} = 0 \ \ \text{for $k \in {\Bbb Z}$}. $$
The doubly infinite Jacobi matrix $\J$  again coincides with $J^-$, which
is now the semi-infinite Jacobi matrix for the Hermite polynomials $H_n(x)$. We
thus find
$$  \lim_{n \to \infty}  \sum_{j=0}^{n-1}
    f\left(\frac{x_j-n^2/2}{n^{3/2}}\right) p_{n-1}^2(x_j;\alpha,n/2,n)
\pi_{j,n} = \int_{-\infty}^\infty f(x) \frac{e^{-x^2}}{\sqrt{\pi}} \, dx,
$$
where now $\pi_{j,n} = \pi_j(\alpha,n/2,n)$,
and in general
$$  \multline
  \lim_{n \to \infty}  \sum_{j=0}^{n-1}
    f\left(\frac{x_j-n^2/2}{n^{3/2}}\right)
    p_{n-k}(x_j;\alpha,n/2,n) p_{n-l}(x_j;\alpha,n/2,n) \pi_{j,n}  \\
 = \frac{1}{\sqrt{2^{k+l-2} (k-1)! (l-1)!}}
   \int_{-\infty}^\infty f(x) H_{k-1}(x)H_{l-1}(x) \frac{e^{-x^2}}{\sqrt{\pi}}
   \, dx .
\endmultline   $$

Weak limits of this kind are useful when studying fluctuation theory
in the Moran (or Bernoulli-Laplace) model in genetics \cite{7}.

\enddocument